\documentclass[12pt]{article}

\usepackage{amsmath, amstext, amsgen, amsbsy, amsopn, amsfonts, amssymb, graphicx,psfrag}
\usepackage[normalem]{ulem}

\usepackage{pdfsync}
\usepackage{epstopdf}
\usepackage{color}

\newtheorem{theorem}{Theorem}[section]

\newtheorem{lemma}[theorem]{Lemma}
\newtheorem{corollary}[theorem]{Corollary}

\title{Involutory Quandles of $(2,2,r)$-Montesinos Links}
\author{
Jim Hoste\\
Pitzer College\\
\\
Patrick D. Shanahan\\
Loyola Marymount University
}
\begin{document}

\maketitle

\begin{abstract}
In this paper we show that  Montesinos links of the form $L(1/2, 1/2, p/q;e)$, which we call $(2,2,r)$-Montesinos links,  have finite involutory quandles.  This generalizes an observation of Winker regarding the $(2, 2, q)$-pretzel links. We  also describe some  properties of these quandles.
\end{abstract}

\section{Introduction}
Let $K$ be a knot or link in $\mathbb S^3$. The knot quandle $Q(K)$ is an algebraic knot invariant that can be derived from a diagram of $K$.  In \cite{JO}, \cite{JO2}, Joyce proves that the knot quandle is a complete invariant of unoriented knots.  However, as with knot groups, it is easy to obtain a presentation of the quandle, but hard to determine if two presentations yield isomorphic quandles.  Passing to the $n$-quandle, $Q_n(K)$, a certain quotient of $Q(K)$, yields a more tractable knot invariant about which much is still unknown. This paper is one piece of a larger project to tabulate all links with finite $n$-quandles and furthermore to explicitly describe those quandles. We hope to complete this project in future publications.

Combining results of Joyce \cite{JO},  \cite{JO2},  and Winker \cite{WI}, it follows that if $Q_n(K)$ is finite, then $\pi_1(\widetilde M_n)$ is finite, where $\widetilde M_n$ is the the $n$-fold cyclic branched cover of $S^3$, branched over $K$. In this paper, we examine the class of Montesinos links $L(1/2, 1/2, p/q;e)$, which we call $(2,2,r)$-Montesinos links, and which are known to have 2-fold cyclic branched covers with finite fundamental group. We show that their 2-quandles, also known as {\it involutory} quandles, are in fact finite. Moreover, we carefully describe the nature of these quandles.
Our main result is the following theorem.

\begin{theorem}\label{main theorem} If $L(1/2, 1/2, p/q;e)$ is a $(2,2,r)$-Montesinos link with $0<p<q$,\\ $\gcd(p,q)=1$, and $e\in \mathbb Z$, then the order of $Q_2(L)$ is $2(q+1)|(e-1)q-p|$.
\end{theorem}

 In Section~\ref{quandles}, we review basic material about involutory quandles and describe the  diagramming method introduced by Winker~\cite{WI} used to find a Cayley graph of a quandle. Then, in Section~\ref{22rMontesinos} we describe the $(2,2,r)$-Montesinos links and derive a presentation for their involutory quandles. We then employ Winker's method in Section~\ref{Cayley graphs} to compute the involutory quandles of these links, describing the quandle by means of its Cayley graph. Finally, in Section~\ref{properties}, we describe several properties of these quandles.

\section{The Involutory Quandle of a Link}\label{quandles}

In this section we quickly review the definition of the fundamental quandle of a link and its associated involutory quandle. We refer the reader to \cite{FR}, \cite{JO}, \cite{JO2}, and \cite{WI} for more detailed information.

A {\it  quandle} is a set $Q$ together with two binary operations $\rhd$ and $\rhd^{-1}$ which satisfies the following three axioms.

\vspace{-.2in}

\begin{itemize}
\item[\bf A1.] $x \rhd x =x$ for all $x \in Q$.
\item[\bf A2.] $(x \rhd y) \rhd^{-1} y = x = (x \rhd^{-1} y) \rhd y$ for all $x, y \in Q$.
\item[\bf A3.] $(x \rhd y) \rhd z = (x \rhd z) \rhd (y \rhd z)$ for all $x,y,z \in Q$.
\end{itemize}

\vspace{-.2in}

 If $L$ is an oriented knot or link in $\mathbb S^3$, then a presentation of its {\it fundamental quandle}, $Q(L)$, can be derived from a regular diagram $D$ of $L$. This process mimics the Wirtinger algorithm. Namely, assign a quandle generator $x_1, x_2, \dots , x_n$ to each arc of $D$, then at each crossing introduce a relation as follows. Suppose the overcrossing strand is labeled $x_j$ while the undercrossing strand to the left (using the orientation of the overcrossing strand) is labeled $x_i$ and the undercrossing strand to the right is labeled $x_k$. Then $x_k \rhd x_j=x_i$. 
 It is easy to check that the three Riedemeister moves do not change the quandle given by this presentation so that the quandle is indeed an  invariant of the oriented link.  

It is important to note that the operation $\rhd$ is, in general, not associative. Throughout this paper, we will assume that any expression without parentheses is to be grouped in left-associated form. For example, the expression $u\rhd^{-1} v\rhd w\rhd x \rhd y\rhd^{-1}  z $ should be taken to be  
$((((u\rhd^{-1} v)\rhd w)\rhd x) \rhd y)\rhd^{-1} z$. It is an interesting exercise to show that the quandle axioms can be used to rewrite every product  in left-associated form. Thus, there is no loss of generality in working entirely with left-associated products.

A quandle is called {\it involutory} if $\rhd = \rhd^{-1}$. In an involutory quandle axiom A2 becomes\\ $(x \rhd y) \rhd y =x$ for all $x, y \in Q$. The involutory quandle for a knot or link is presented in the same manner as its fundamental quandle except that the orientation of the link is now irrelevant.  The involutory quandle of a link $L$ will be denoted by $Q_2(L)$. The involutory quandle of a link may be thought of as a quotient of the fundamental quandle and in general one should expect a significant loss of information in passing from $Q(L)$ to $Q_2(L)$. All quandles  considered in this paper will be involutory quandles.
Finally, we adopt the exponential notation used in \cite{FR}, denoting $x \rhd y$ as $x^y$. This has the advantage of further clarifying the ambiguity caused by lack of associativity. Now $x^{yz}$ will be taken to mean $(x^y)^z=(x \rhd y)\rhd z$ whereas $x^{y^z}$ will mean $x\rhd (y \rhd z)$.

The following lemma from \cite{FR}, which  describes how to re-associate a product in an involutory quandle, will be used repeatedly in this paper. 
\begin{lemma} If $x, y, z_1, z_2, \dots, z_k$ are elements of an involutory quandle, then
$$x^{y^{z_1 z_2 \dots z_k}} = x^{z_k \dots z_2 z_1 y z_1 z_2 \dots z_k}.$$
\label{leftassoc}
\end{lemma}
\vspace{-.5 in}
One of the main tools we use in this paper is Winker's  method to produce the Cayley graph of an  involutory quandle from a given presentation \cite{WI}. The method is similar to the well-known Todd-Coxeter method of enumerating cosets of a subgroup of a group \cite{TC}. Suppose $Q_2(L)$ is presented as
$$Q_2=\langle x_1, x_2, \dots, x_n \, |\, x_{j_1}^{w_1}=x_{k_1}, \dots, x_{j_m}^{w_m}=x_{k_m} \rangle,$$
  where each $w_i$ is a word in the generators. If $y$ is any element of the quandle, then it follows from Lemma~\ref{leftassoc} that $$y^{\overline{w_i}x_{j_i}w_i}=y^{x_{k_i}},$$ where $\overline{w}$ is the word $w$ written backwards. Winker calls this relation the {\it secondary relation}  associated to the {\it primary relation} $x_{j_i}^{w_i}=x_{k_i}$. Winker's method now proceeds as follows.
  
  \begin{enumerate}
  \item Begin with one vertex for each generator $x_i$, labeled with the integer $i$. 
  \item Add a loop at each vertex $x_i$,  labeled $x_i$. (This encodes Axiom A1.)
  \item Using each primary relation $x_{j_i}^{w_i}=x_{k_i}$,  in their given order, introduce vertices and edges as necessary to create the path from $x_{j_i}$ to $x_{k_i}$ given by $w_i$. This process is called {\it tracing} the primary relation. Vertices are labeled with consecutive integers in the order they are introduced and edges are labeled with the corresponding generator. Tracing the $i$-th relation may introduce like-labeled edges sharing a vertex which we now  identify, possibly leading to other necessary identifications. This process is called {\it collapsing} and all collapsing is carried out before tracing the $(i+1)$-st relation. 
\item Proceeding in order through the list of vertices, trace and collapse each secondary relation (in order). All secondary relations are traced and collapsed at vertex $i$ before proceeding to vertex $i+1$.
 \end{enumerate}
  
If the method terminates in a finite graph, then this is the Cayley graph of the quandle. However, the method may not terminate. The reader is referred to Winker~\cite{WI} for more details. 

\begin{figure}[htbp]
\vspace*{13pt}
\centerline{\includegraphics*[scale=.6]{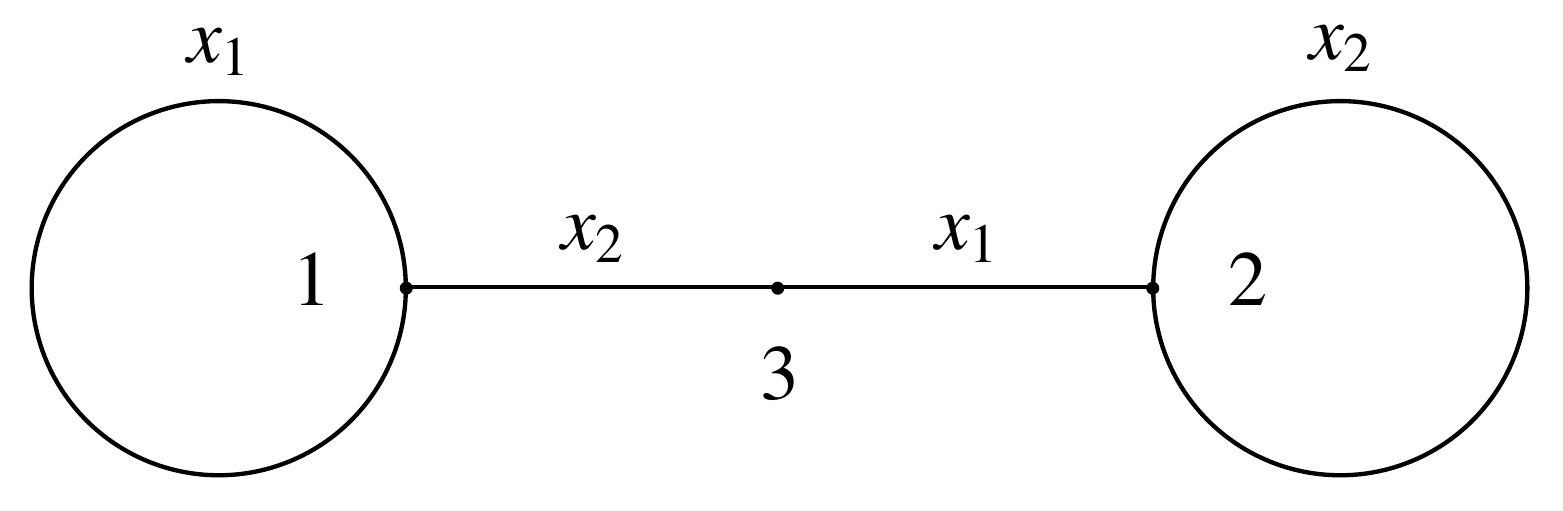}}
\caption{The Cayley graph of the involutory quandle of the right-handed trefoil knot.}
\label{trefoil}
\end{figure}

For example, let $K$ be the the standard 3-crossing diagram of the right-handed trefoil knot. If the arcs are labeled $x_1, x_2$, and $x_3$ in order as we traverse the knot, then we obtain the presentation
$$\langle x_1, x_2, x_3 \, |\, x_1^{x_3}=x_2, x_2^{x_1}=x_3, x_3^{x_2}=x_1\rangle$$
of its involutory quandle.
As explained in~\cite{FR}, presentations can be changed by the analog of Tietze transformations. Here we can eliminate $x_3$ obtaining  
$$\langle x_1, x_2 \, |\, x_1^{x_2}=x_2^{x_1}\rangle,$$
to which we now apply Winker's method.
 After the first two steps, we arrive at the graph shown in Figure~\ref{trefoil}. Tracing the single secondary relation, $q^{x_1x_2x_1x_2x_1x_2}=q$ at each of the vertices (in order) introduces no additional vertices and requires no collapsing. Therefore, this is the Cayley graph of the involutory quandle of the right-handed trefoil knot.

\section{The  $(2,2,r)$-Montesinos Links}\label{22rMontesinos}

\begin{figure}[h]
\vspace*{13pt}
\centerline{\includegraphics[scale=.55]{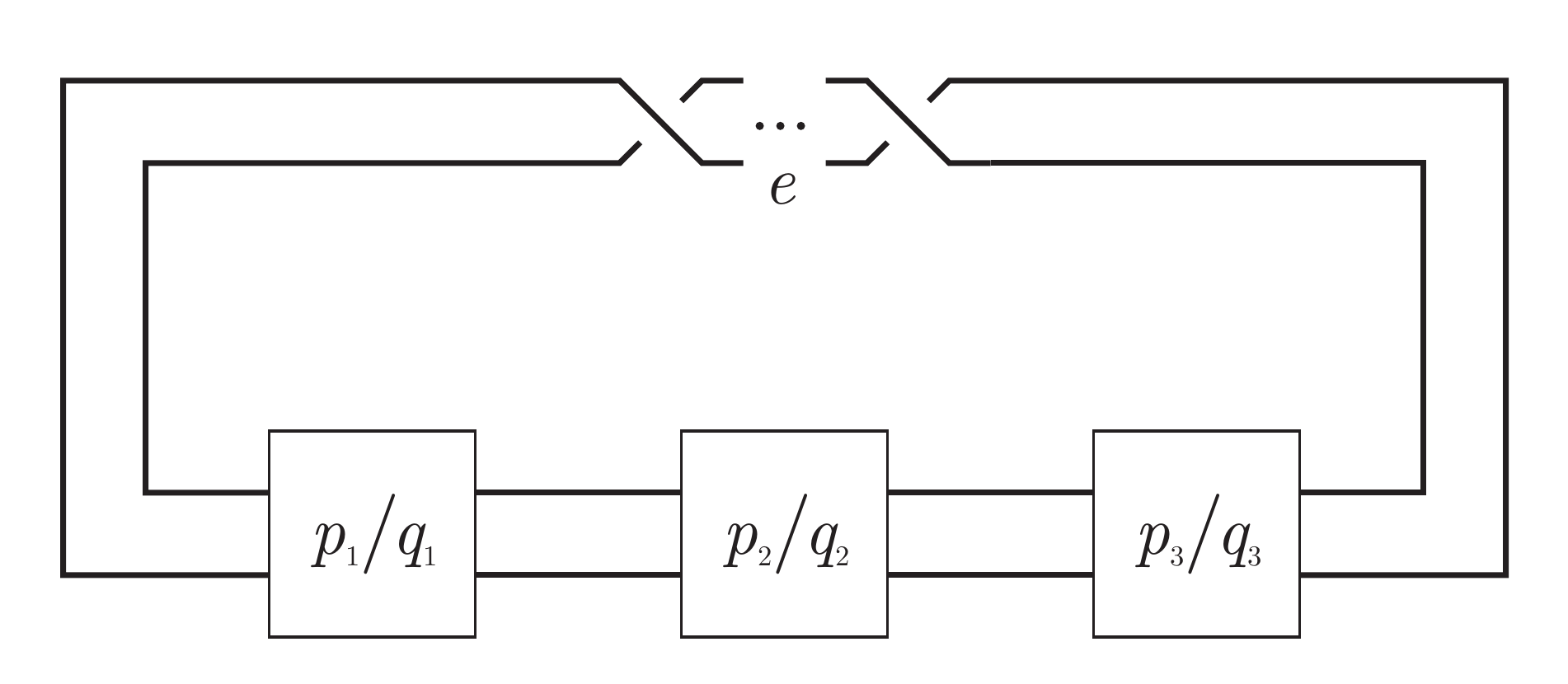}}
\caption{The Montesinos link $L(p_1/q_1, p_2/q_2, p_3/q_3; e)$}
\label{generalMontesinos}
\end{figure}

Consider the unoriented Montesinos link $L(p_1/q_1, p_2/q_2, p_3/q_3; e)$, shown in Figure~\ref{generalMontesinos}, where $q_i>0$, $\gcd(p_i,q_i)=1$, and $e \in \mathbb Z$.  Here the box labelled $p_i/q_i$ denotes the rational tangle determined by the fraction $p_i/q_i$ and $e$ is the number of right-handed half-twists.  By making use of flypes, we may change each numerator $p_i$ by multiples of $q_i$ at the expense of changing $e$. In particular, if some $p_i$ is replaced with $p_i \pm q_i$, then we replace $e$ with $e \pm 1$. Hence by allowing arbitrary values for $e$, we may assume that $0<p_i<q_i$ for each tangle. In addition, cyclic permutation of the tangles together with rotation around a vertical axis allows us to assume $0<q_1 \le q_2 \le q_3$. The reader is directed to \cite{BZ} for more information on Montesinos links.
For the remainder of the paper we consider the link  $L=L(1/2,1/2,p/q;e)$. We begin by deriving a presentation for the involutory quandle $Q_2(L)$. 

\begin{figure}[h]
\vspace*{13pt}
\centerline{\includegraphics[scale=.5]{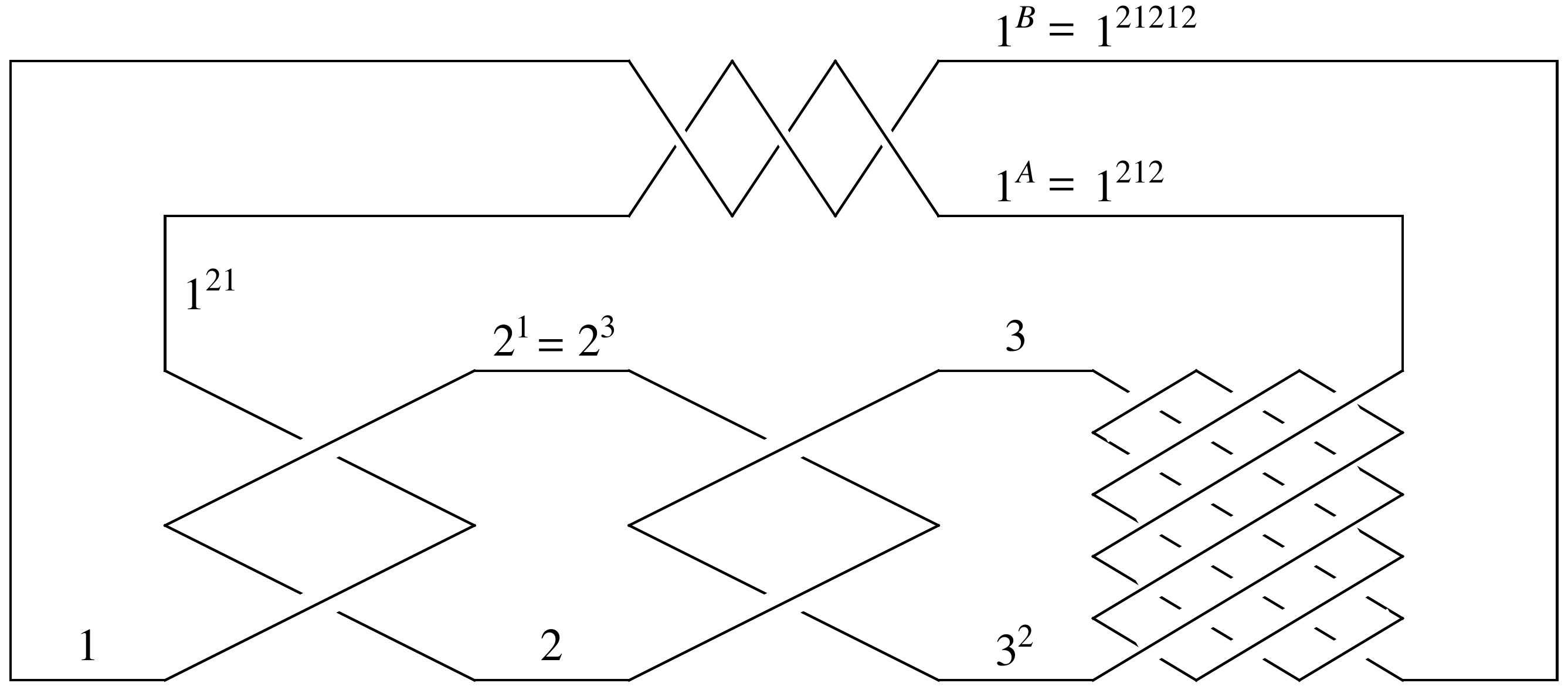}}
\caption{The Montesinos link $L(1/2, 1/2, 3/5; e)$}
\label{montesinosLink}
\end{figure}

We illustrate this process with the link $L(1/2, 1/2, 3/5; 3)$ shown in Figure~\ref{montesinosLink}.  Label the arcs of the diagram with quandle elements 1, 2, and 3 as shown in the figure. (These might more properly be called $x_1, x_2$, and $x_3$, but we simply write 1, 2, 3.) The elements 1, 2, and 3 generate the quandle. Notice that for all $p/q$, we obtain the relation $2^1=2^3$. However, there are two more relations that depend on both $e$ and $p/q$. First note that by moving from left to right through the $e$ right-handed half-twists in the upper part of the diagram, the two arcs to the right receive labels of $1^A$ and $1^B$ where $A$ and $B$ are some words in the generators 1 and 2 that depend on $e$. Next, if we follow the strand labeled $3^2$ into the $p/q$ tangle, we will come out at one of the other four endpoints of the tangle and obtain a second relation. Similarly, following the other string in this 2-string tangle gives the third relation.

We begin by determining $A$ and $B$, which are easily found by induction on $e$. In the following lemma, and throughout the paper, if $X$ is a word in the generators, then we use the notation $\overline X$ to denote $X$ written backwards. Additionally, we use $[X]^n$   to stand for $X$ repeated $n$ times when $n>0$,  to stand for the empty string when $n=0$, and to stand for $[\overline X]^{-n}$ when $n<0$.
\begin{lemma}The words $A$ and $B$ depend on $e$ and are given as, $A=\overline A=[21]^{e-2}2$ and $B=\overline B=[21]^{e-1}2$.
\label{AB}
\end{lemma}

We now consider the $p/q$ tangle. There are three different cases: $q$ even and  $p$ odd, $q$ odd and $p$ even, and finally, $q$ odd and $p$ odd. The Montesinos link has a total of three components if $q$ is even and two if $q$ is odd. Moreover, the parity of $p$ and $q$ determine which pairs of endpoints are connected by the two strings of the tangle.  If $q$ is even, then the lower left endpoint is connected to the upper left endpoint. If $q$ is odd, and $p$ is even, then the lower left endpoint is connected to the lower right endpoint. Finally, if both are odd, the lower left endpoint is connected to the upper right endpoint.

\begin{figure}[h ]
\vspace*{13pt}
\centerline{\includegraphics[scale=.5]{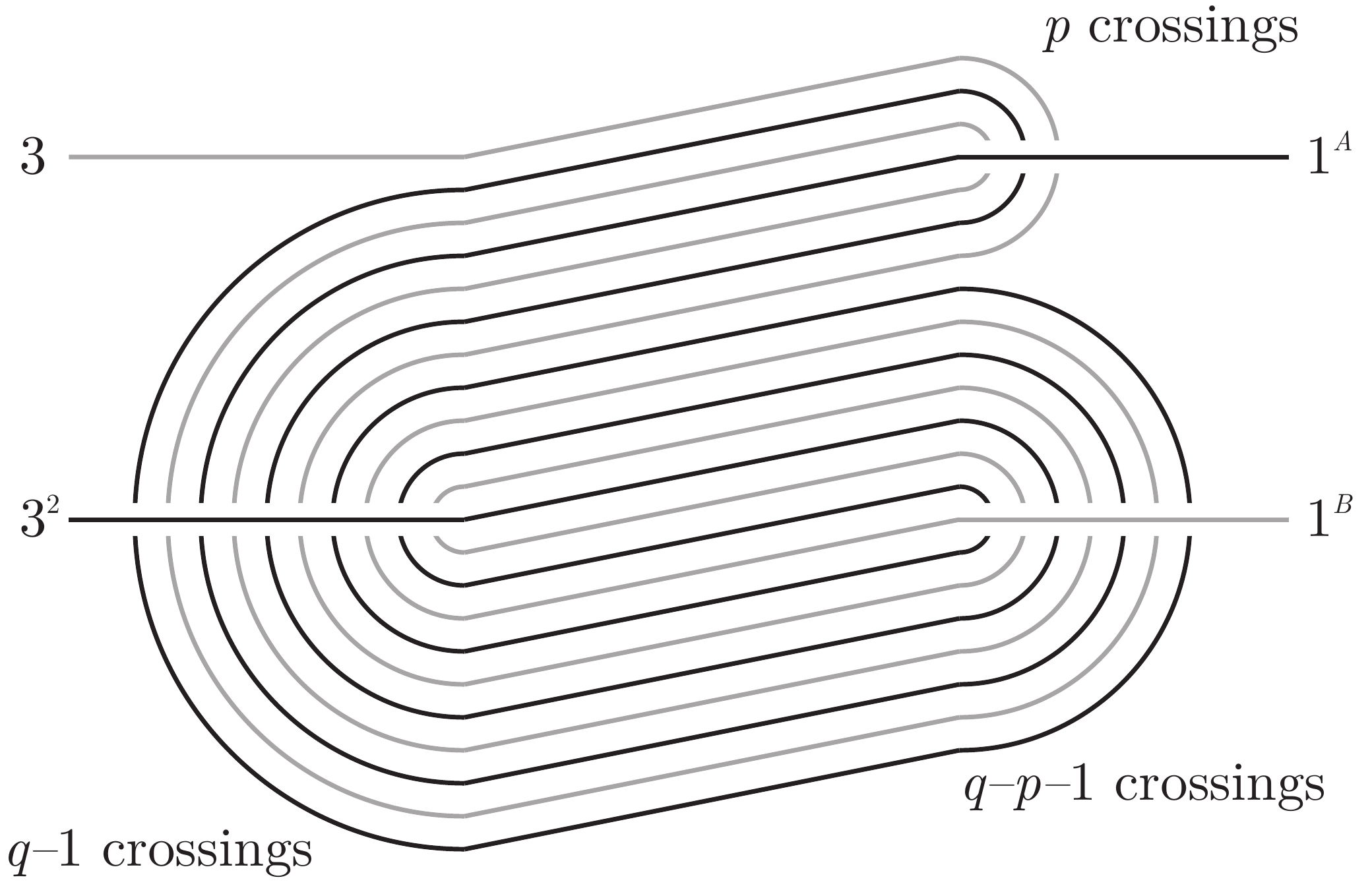}}
\caption{The case $p$ and $q$ both odd, shown with $p/q=3/11$.}
\label{pqtangle}
\end{figure}

We describe the case with $p$ and $q$ both odd.   After isotopy, the tangle can be represented as shown in Figure~\ref{pqtangle}.  Consider the lower left endpoint of the tangle. This arc is labeled $3^2$. As we move along this arc, we alternately pass under arcs labelled $1^A$ or $1^B$ and $3^2$. We pass under the arc labelled $3^2$ a total of $(q-1)/2$ times, under $1^A$ a total of $(p-1)/2$ times, and under $1^B$ a total of $(q-p)/2$ times. This path ends at $1^A$ yielding a relation of the form $3^{2X} = 1^A$. By Lemmas~\ref{leftassoc} and \ref{AB}, $X$ is an alternating product of the words $A1A$ or $B1B$ and $232$. For example, with $p/q=3/11$ as shown in Figure~\ref{pqtangle} we obtain the relation
$$3^{2B1B232A1A232B1B232B1B232B1B232}=1^A.$$
For any $p/q$, we can simplify the relation using the following result, the proof of which we postpone to the next section. 
\begin{lemma}\label{simplify relation} For all $x \in Q_2(L)$, $x^{A1A232B1B} = x^{B1B232A1A}$.
\end{lemma}

Repeated application of this lemma allows us to collect all $A1A232$ and $B1B232$ expressions together and so the relation $3^{2X}=1^A$ can be rewritten as 
$$3^{2[A1A232]^{(p-1)/2}[B1B232]^{(q-p)/2}}=1^A.$$
Now applying this same process starting at the arc labelled $3$ gives a second relation 
$$3^{[A1A232]^{(p+1)/2}[B1B232]^{(q-p-2)/2}}=1^B.$$
Finally, a similar analysis for the other parities of $p$ and $q$ gives the following presentation for the involutory quandle. The reader is warned that, in an effort to make the relations appear more uniform across the cases, we introduce some canceling generators to rewrite the last case.
 
\begin{lemma}\label{R2R3} If $L$ is the Montesinos link $L(1/2, 1/2, p/q; e)$ with $0<p<q$, then the involutory quandle of $L$ can be presented as 
$Q_2(L)=\langle 1, 2, 3\, | \, R_1,R_2, R_3\rangle$, where  $R_1$ is the relation $2^1=2^3$.  The relations $R_2$ and $R_3$ depend on the parity of $p$ and $q$ as follows:
\vskip -.35 in
\renewcommand*{\arraystretch}{1.5}
$$\begin{array}{llll}
p&q&R_2&R_3\\
\hline
\text{odd}&\text{odd}&
3=1^{A[232 B1B]^{(q-p)/2}[232A1A]^{(p-1)/2}2}&
3=1^{B [232 B1 B]^{(q-p-2)/2}[232A1A]^{(p+1)/2}}\\

\text{even}&\text{odd}&
3=1^{B[232B1B]^{(q-p-1)/2}[232A1A]^{p/2}2}&
3=1^{A[232B1B]^{(q-p-1)/2}[232A1A]^{p/2}}\\

\text{odd}&\text{even}&
3=3^{2 [232B1B]^{(q-p-1)/2} [232A1A]^{(p+1)/2} }&
1=1^{A[232B1B]^{(q-p+1)/2}[232 A1A]^{(p-1)/2}B}\\
\end{array}
$$
\end{lemma}

\section{The Cayley Graph of $Q_2(L)$}\label{Cayley graphs}

In this section we explicitly construct the Cayley graph of $Q_2(L(1/2, 1/2, p/q; e))$ with respect to  the generators 1, 2, and 3, given in Lemma~\ref{R2R3}. Let ${\cal C}_i$ be the component of the Cayley graph containing the generator $i$.  Each component ${\cal C}_i$ is best described as the quotient of an infinite, labeled, planar graph $ {\cal T}_i$ which we now describe. Start with a single rectangular tile one unit wide and four units tall with its lower left corner located at the point $(0,1)$. Consider this rectangle as a graph with six edges and six vertices located at $(0,1)$, $(0,3)$, $(0,5)$, $(1,5)$, $(1,3)$, $(1,1)$.  Let $ {\cal T} $ be the union of all translates of this tile obtained by shifting $m$ units horizontally and $n$ units vertically, where $m$ and $n$ are integers with $n \equiv 2m \text{ (mod 4)}$. Label all horizontal edges with 2, and all vertical edges whose upper endpoint is located at $(x,y)$ with 3 if $y \equiv 3$ (mod 4) and with 1 if $y \equiv 1$ (mod 4). Finally, let ${\cal T}_i$ be obtained from $\cal T$ by labeling the vertex $(0,1)$ with generator $i$.

The graph ${\cal T}$ admits a number of label-preserving symmetries that we will use to describe $Q_2(L)$. For any integers $n$ and $m$, let $\rho_{(n, 2m)}$ denote the $\pi$-radian rotation centered at $(n,2m)$, let $\tau_{(2n, 4m)}$ denote the translation $(x,y) \to (x+2n, y+4m)$, let $\phi_{n+1/2}$ denote reflection across the line $x=n+1/2$, and let $\gamma_{(4n,m+1/2)}$ denote the glide reflection $\tau_{(0,4n)}\circ \phi_{m+1/2}$. It is easy to check that each of these transformations is a label-preserving automorphism of $\cal T$. If $\psi_1, \psi_2, \dots, \psi_k$ are label-preserving automorphisms of $\cal T$ we denote by $\langle \psi_1, \psi_2, \dots, \psi_k \rangle$ the group of automorphism which they generate and by ${\cal T}/\langle \psi_1, \psi_2, \dots, \psi_k \rangle$ the quotient graph.
    
\begin{theorem}\label{Cayley graph} Suppose $L=L(1/2, 1/2, p/q;e)$ and let $w=(e-1)q-p$. Then the components of the Cayley graph of $Q_2(L)$ are:
\vspace{-.25 in}
\begin{enumerate}
\item[(a)] If $q$ is odd, then ${\cal C}_1={\cal C}_3 ={\cal T}_1/\langle \rho_{(0,0)},\rho_{(w,0)},\rho_{(0,2q)} \rangle$ and ${\cal C}_2 = {\cal T}_2/\langle \tau_{(0,4)}, \phi_{1/2}, \phi_{w+1/2} \rangle$. 
\item[(b)]If $q$ is even, then ${\cal C}_1 ={\cal T}_1/\langle \rho_{(0,0)},\rho_{(0,2q)}, \gamma_{(2q, w/2)} \rangle$, ${\cal C}_2 = {\cal T}_2/\langle \tau_{(0,4)}, \phi_{1/2}, \phi_{w+1/2} \rangle$, and ${\cal C}_3={\cal T}_3/\langle \rho_{(0,2)},\rho_{(0,2+2q)}, \gamma_{(2q, w/2)} \rangle$.  The component ${\cal C}_3$ is obtained by trading all labels of 1 and 3 in ${\cal C}_1$.
\end{enumerate}
\end{theorem}

{\bf Remark:} It follows from Theorem~\ref{Cayley graph} that a fundamental domain ${\cal F}_1$ for ${\cal C}_1$ can be taken to be the subset of ${\cal T}_1$ contained in $[0, |w|]\times[0,4q]$ if $q$ is odd and contained in $[-|w|/2, |w|/2] \times [0, 2q]$ if $q$ is even. If $q$ is odd, ${\cal F}_1$ is not quite a fundamental domain, as pairs of vertices symmetric with respect to $(0,2q)$ or with respect to $(|w|,2q)$ still need to be identified.  A fundamental domain ${\cal F}_2$ for  ${\cal C}_2$, in either case, is the subset of ${\cal T}_2$ contained in $[1/2, |w|+1/2]\times[0,4]$. Counting the number of vertices in ${\cal F}_1$ and ${\cal F}_2$ gives the following corollary.

\begin{corollary} \label{corollary}The orders of the components of the Cayley graph of $Q_2(L)$ are:
\vspace{-.25 in}
\begin{enumerate}
\item[(a)] If $q$ is odd, then $|{\cal C}_1| =2 q |w|$ and $|{\cal C}_2|=2|w|$. 
\item[(b)]If $q$ is even, then $|{\cal C}_1|=|{\cal C}_3| = q |w|$ and $|{\cal C}_2|=2|w|$. 
\end{enumerate}
\vspace{-.25 in}
In both cases, $|Q_2(L)|=2(q+1)|w|$.
\end{corollary}
Notice that Theorem~\ref{main theorem} now follows immediately from Theorem~\ref{Cayley graph} and its corollary.

\begin{figure}[htbp]
\vspace*{13pt}
\centerline{\includegraphics*[scale=.9]{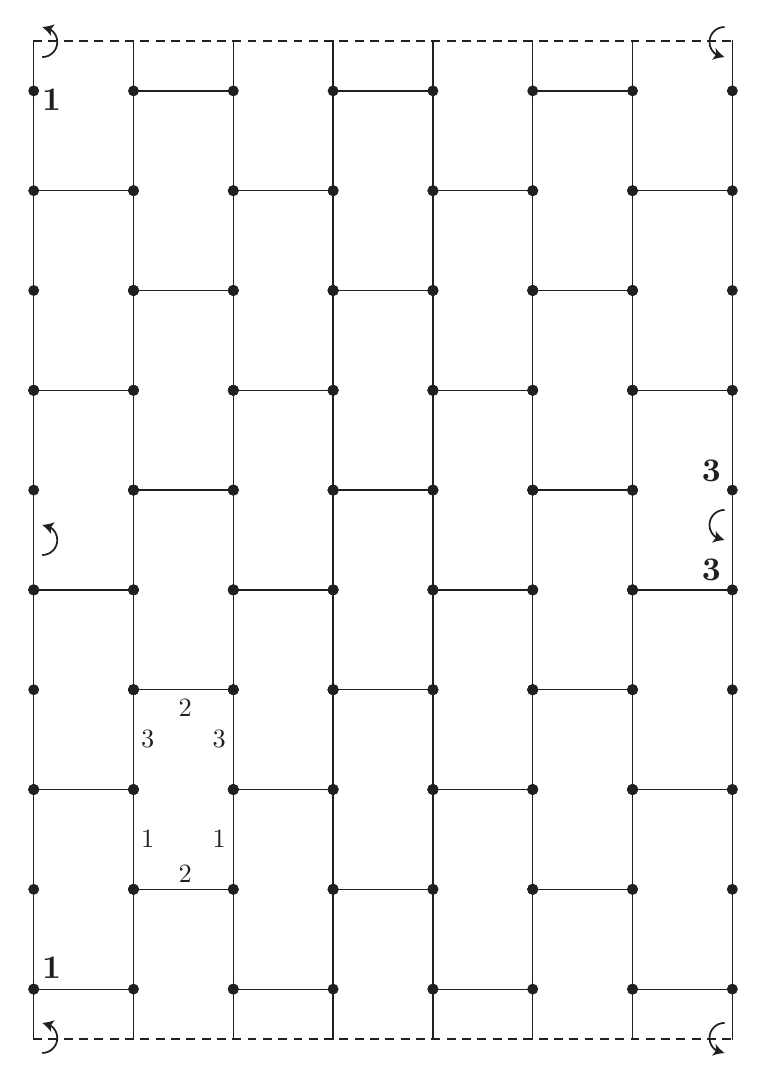}}
\caption{The fundamental domain ${\cal F}_1$ of $Q_2(L)$ for $p=3$, $q=5$, and $e=3$.}
\label{fund domain}
\end{figure}

\begin{figure}[htbp]
\vspace*{13pt}
\centerline{\includegraphics*[scale=.5]{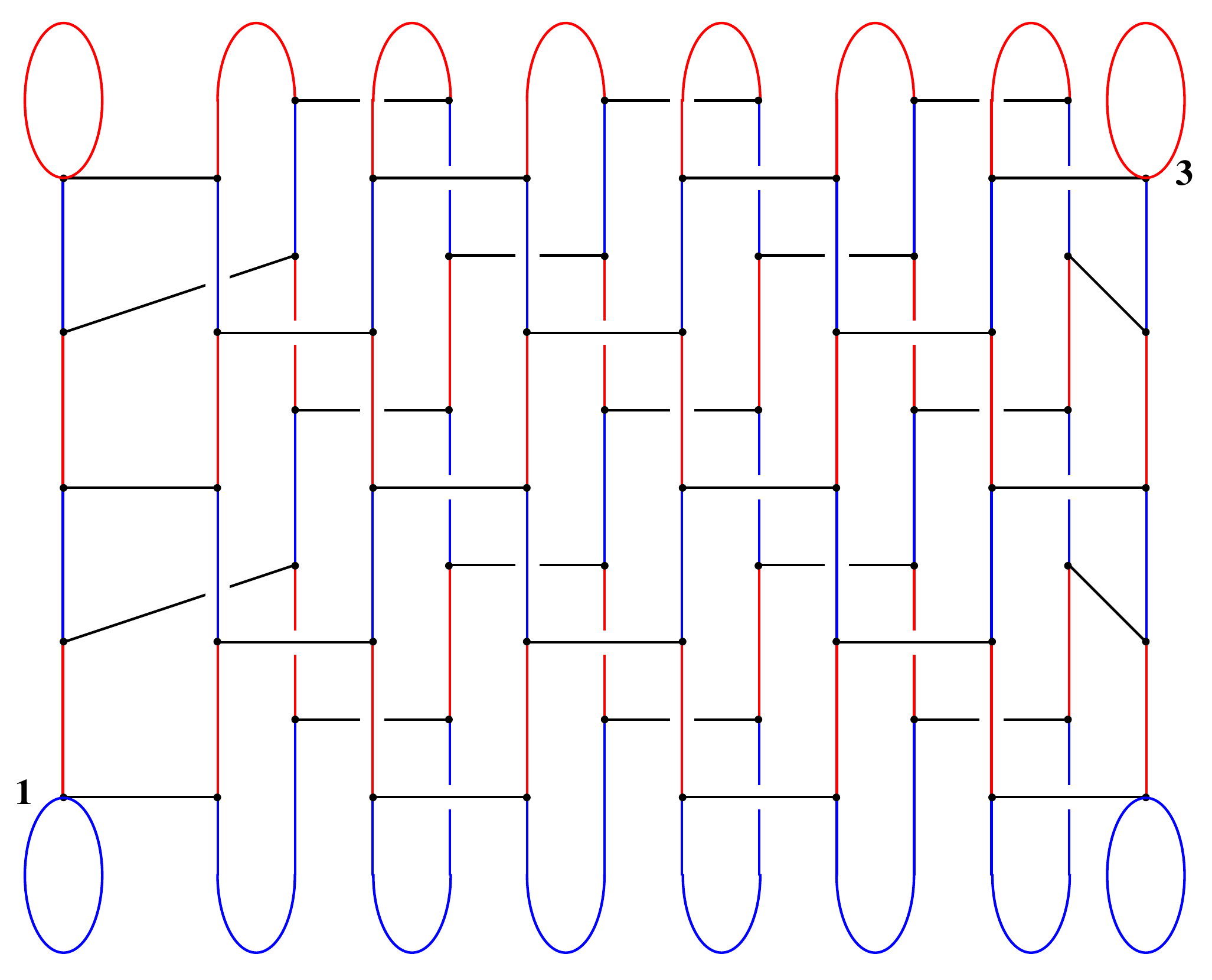}}
\caption{The component ${\cal C}_1$ of $Q_2(L)$ for $p=3$, $q=5$, and $e=3$. Here blue, black, and red correspond to the generators 1, 2, and 3, respectively.}
\label{C1 q odd}
\end{figure}

\begin{figure}[htbp]
\vspace{0 in}
\centerline{\includegraphics*[scale=.5]{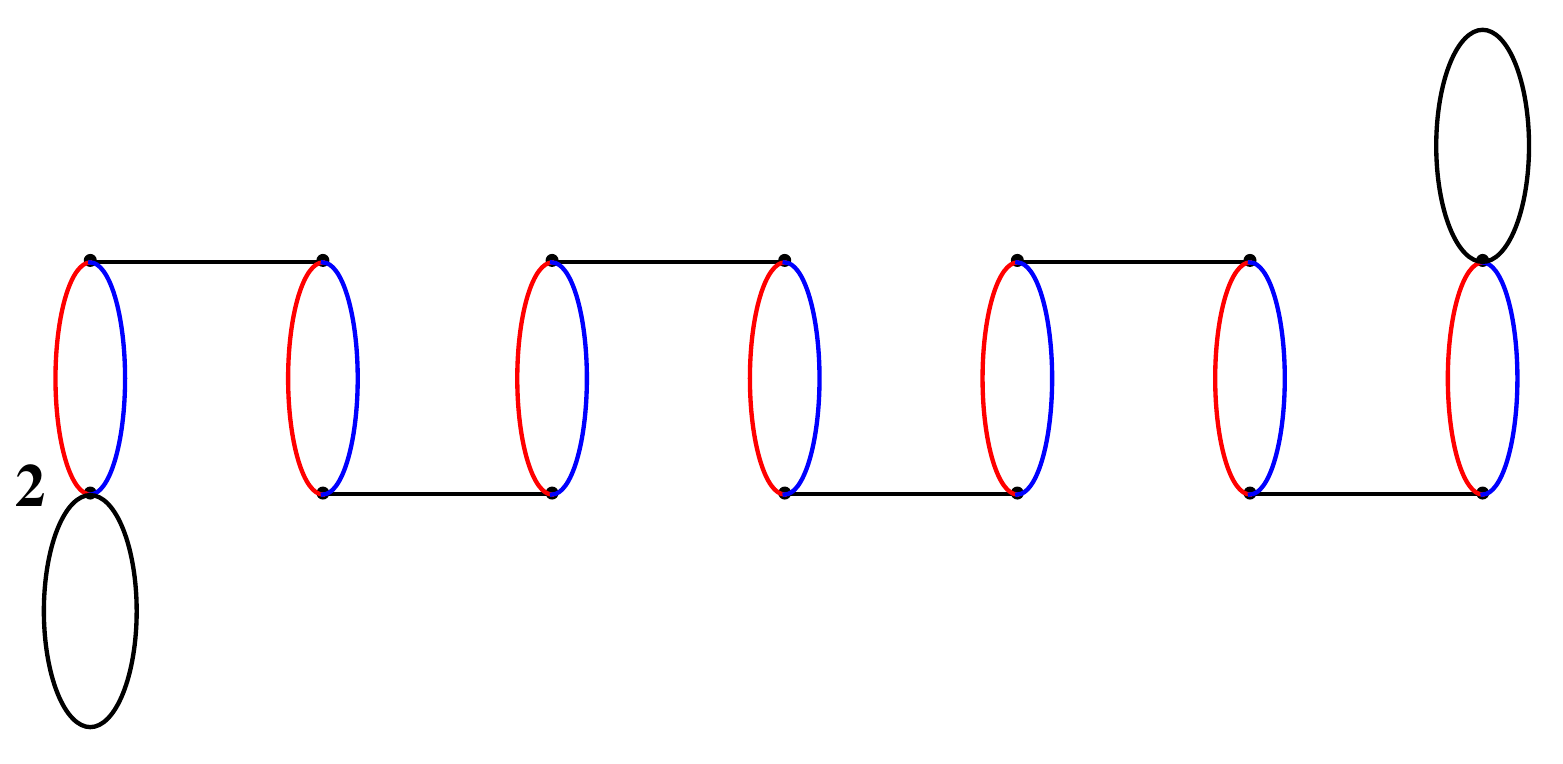}}
\vspace{0 in}
\caption{The component ${\cal C}_2$ of $Q_2(L)$ for $p=3$, $q=5$, and $e=3$. Here blue, black, and red correspond to the generators 1, 2, and 3, respectively.}
\label{C2}
\end{figure}

\begin{figure}[htbp]
\vspace{-.5 in}
\centerline{\includegraphics*[scale=.4]{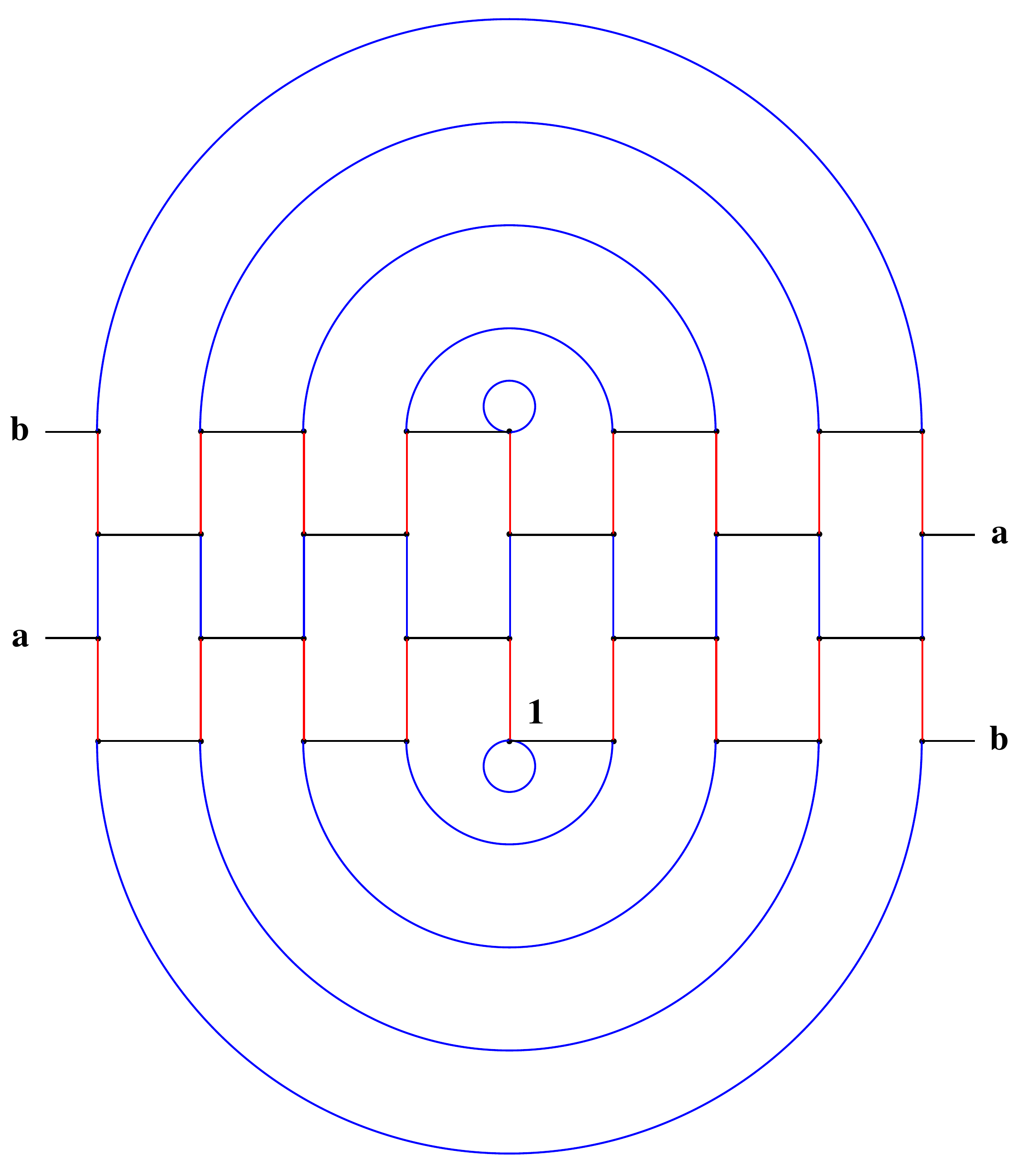}}
\vspace{0 in}
\caption{The component ${\cal C}_1$ of $Q_2(L)$ for $p=3$, $q=4$, and $e=4$. Here blue, black, and red correspond to the generators 1, 2, and 3, respectively. Letters indicate how half-edges at the right and left join. }
\label{C1 q even}
\end{figure}

For example, in the case where $p/q=3/5$ and $e=3$, ${\cal F}_1$, ${\cal C}_1$,  and ${\cal C}_2$ are shown in Figures~\ref{fund domain}, \ref{C1 q odd}, and \ref{C2}. In Figure~\ref{C1 q even}, we illustrate the case of $p/q=3/4$ with $e=4$. 

Our overall strategy for proving these theorems is to apply the diagramming method of Winker~\cite{WI} to a  presentation of the quandle $Q_2(L)$.   In order to simplify the application of Winker's method, we will first pass to an equivalent presentation of the quandle. Our approach is motivated by first considering the involutory quandle $\widetilde Q_2=\langle 1, 2, 3\,|\, 2^1=2^3\rangle$. Imposing additional relations on $\widetilde Q_2$ cause vertices in its Cayley graph to be identified, thus leading to the Cayley graph of $Q_2(L)$.
 
Applying Winker's method to $\widetilde Q_2$, the graph $ {\cal T}$ arises in a natural way because the primary relation $2^1=2^3$ gives rise to the secondary relation $x^{121}=x^{323}$ for every element $x$ in $\widetilde Q_2$. Thus it is not surprising to see the hexagons of $ {\cal T}$ appearing in the end result. In fact, Winker's method results in  three  components: 
\begin{align*}
\widetilde {\cal C}_1&={\cal T}_1/\langle \rho_{(0,0)}\rangle\\
\widetilde {\cal C}_2&={\cal T}_2/\langle \tau_{(0,4)}, \phi_{1/2}\rangle\\
\widetilde {\cal C}_3&={\cal T}_3/\langle \rho_{(0,2)} \rangle
\end{align*}

Starting with the graphs ${\cal T}_1$ and ${\cal T}_2$,  and applying the relations $R_2$ and $R_3$ in Lemma~\ref{R2R3}, we were led  to Theorem~\ref{Cayley graph}. However, we will now treat this approach as simply providing guidance to our formal proof of the theorem. Using $ {\cal T}$, we will see how to rewrite the presentation of $Q_2(L)$ in such a way that applying Winker's method to the new presentation will produce the Cayley graphs described in Theorem~\ref{Cayley graph} after tracing and collapsing the primary relations. Each of the secondary relations will then hold without further modification of the graph. 

We begin with some general observations about tracing a word in the generators starting at various vertices in $ {\cal T}$. First notice that there are four {\it types} of vertices $(x,y)$ in ${\cal T}$ determined by the residues $x$ (mod 2) and $y$ (mod 4). We denote these four vertex types in the obvious way as {\bf 01}, {\bf 03}, {\bf 11}, and {\bf 13}.  

Given a word $X$ in the generators 1, 2, and 3, we can start at any vertex $a$ in $ {\cal T}$, trace the word $X$, and arrive at a vertex $b$.  It is easy to see that the vertex type of $b$ depends only on $X$ and the type of $a$. Thus, each word $X$ determines a permutation of the set of four types of vertices. Let $\sigma(X)$ denote this permutation. For example, $\sigma(2)$ is the product of two transpositions, written in cycle notation as  $\sigma(2)=({\bf 01}\, {\bf 11})({\bf 03}\, {\bf 13})$.

If tracing the word $X$ starting at vertex $a$ in $ {\cal T}$ ends in vertex $b$, then let $\Delta_a X$ denote the  vector $b-a$ in the plane. We call $\Delta_a X$ the {\it displacement} of $X$ starting at $a$.  Because of the translational symmetry taking any vertex to another one of the same type, it follows that $\Delta_a X$ depends only on $X$ and the type of $a$. Therefore, we will often write one of {\bf 01}, {\bf 03}, {\bf 11}, or {\bf 13} in place of $a$. The following Lemma will be helpful in determining displacements. The proof is left to the reader.

\begin{lemma}
\label{displacements} 
Let $X$ be any word in the generators $1$, $2$, and $3$. If $\Delta_{\bf 01} X= \langle x, y \rangle$, then
\vspace{-.2 in}
\begin{enumerate}
\item[(a)] 
$\Delta_{\bf 03}X=\langle -x, -y \rangle,\,\,
\Delta_{\bf 11}X=\langle -x, y \rangle,\,\,\text{and }
\Delta_{\bf 13}X=\langle x, -y \rangle.$
\item[(b)] $\Delta_{\bf 01}\overline{X}$ equals $\langle -x, -y \rangle$, $\langle x, y \rangle$, $\langle x, -y \rangle$, or $\langle -x, y \rangle$, if $\sigma(X)$ equals $\text{id}$, $\sigma(1)$, $\sigma(2)$, or $\sigma(12)$, respectively.
\end{enumerate}
\end{lemma}

Lemma~\ref{displacements} implies that if the words $X$ and $Y$ have the same displacement starting at some vertex $a$, then they have the same displacement starting at any vertex. If this is the case, we say that $X$ and $Y$ {\it have the same displacement}.
Using Lemma~\ref{displacements}, we may compute the permutations and displacements of the words that appear in the relations $R_2$ and $R_3$. Again, details are left to the reader.

\begin{lemma} \label{specific displacements} The permutation and displacement of $A$, $B$, $232A1A$, and $232B1B$ are:
\begin{align*}
\sigma(A)&=\left \{ \begin{array}{ll}\sigma(2), &\text{ if $e$ is even}\\
\sigma(1), &\text{ if $e$ is odd} \end{array} \right. & \Delta_{\bf 01}(A)&=
\left \{ \begin{array}{ll}\langle e-1, 0\rangle, &\text{ if $e$ is even}\\
\langle e-1, -2\rangle, &\text{ if $e$ is odd}\end{array} \right.\\
\sigma(B)&=\left \{ \begin{array}{ll}\sigma(1), &\text{ if $e$ is even}\\
\sigma(2), &\text{ if $e$ is odd} \end{array} \right. & \Delta_{\bf 01}(B)&=
\left \{ \begin{array}{ll}\langle e,-2\rangle, &\text{ if $e$ is even}\\
\langle e,0\rangle, &\text{ if $e$ is odd}\end{array} \right.\\
\sigma(232A1A)&=\text{id} & \Delta_{\bf 01}(232A1A)&=\langle 4-2e, 4 \rangle\\
\sigma(232B1B)&=\text{id} & \Delta_{\bf 01}(232A1A)&=\langle 2-2e, 4 \rangle\\
\end{align*}
\end{lemma}
\vspace{-.5 in}

Let us now focus on the case where both $p$ and $q$ are odd. The component  $\widetilde{\cal C}_1$  is the quotient  ${\cal T}_1/\langle \rho_{(0,0)}\rangle$. The pair of vertices $(0,1)$ and $(0,-1)$ both represent the generator 1. The relation $R_2$ is of the form $1^\alpha=3$ where $\alpha=A[232 B1B]^{(q-p)/2}[232A1A]^{(p-1)/2}2$.  Thus,  if we start at the vertex $(0,1)$ and trace the word $\alpha$, we will arrive at a vertex $b$ that is equal to the generator 3. To determine the coordinates of the vertex $b$, we use Lemmas~\ref{displacements} and \ref{specific displacements} to calculate the displacement of $\alpha$ starting at $(0,1)$. Assume that $e$ is even. Then, by Lemma~\ref{specific displacements}, we have $\sigma(A)=\sigma(2)$ and $\Delta_{\bf 01}A = \langle e-1,0 \rangle$. So, starting at $(0,1)$ and following $A$,  we end at vertex $(e-1, 1)$ which has type $\bf 11$. This observation together with the fact that $\sigma(232A1A)$ and $\sigma(232B1B)$ are the identity then gives
\begin{align*}
\Delta_{\bf 01}\alpha&=\Delta_{\bf 01} A+\frac{q-p}{2} \, \Delta_{\bf 11}232B1B+ \frac{p-1}{2} \, \Delta_{\bf 11}232A1A+\Delta_{\bf 11}2
\\
\\
&=\langle e-1, 0 \rangle + \frac{q-p}{2} \, \langle 2e-2, 4 \rangle + \frac{p-1}{2} \, \langle 2e-4,4 \rangle + \langle -1,0 \rangle\\
\\
&=\langle (e-1)q-p, 2q-2 \rangle\\
&=\langle w, 2q-2 \rangle. 
\end{align*}
Thus, the vertex $(w, 2q-1)$ is equal to the generator 3. Since $p$ and $q$ are odd and $e$ is even this point is of type {\bf 01}. Moving from this vertex along the edge labelled $3$ we arrive at the type $\bf 03$ vertex $(w, 2q+1)$.  However, because $3^3=3$, it follows that the vertex $(w, 2q+1)$ also represents the generator 3.  If $e$ is odd, then a similar calculation shows that $\Delta_{\bf 01}\alpha=\langle w, -2q\rangle$ and that the vertices $(w, -2q\pm1)$ represent the generator 3.

The relation $R_3$ is of the form $1^\beta=3$ where $\beta=B [232 B1 B]^{(q-p-2)/2}[232A1A]^{(p+1)/2}$. Calculating the displacement of $\beta$ yields $\Delta_{\bf 01}\beta=\langle w, -2q\rangle$ if $e$ is even and $\Delta_{\bf 01}\beta=\langle w, 2q-2\rangle$ if $e$ is odd. Therefore, regardless of the parity of $e$, we have that the vertices $(w, 2q\pm1)$ and $(w, -2q\pm1)$ all represent the generator 3.
It now follows that any two vertices that are symmetric with respect to either $(w,  2q)$ or $(w, -2q)$ represent equivalent elements in $Q_2(L)$. Thus these points become rotational centers of symmetry.

Similar calculations give the displacements of the words $\alpha$ and $\beta$ that appear in relations $R_2$ and $R_3$ for the other parity cases. We summarize the results in Table~\ref{R2R3 displacements}.

\begin{table}[t]
\caption{Displacements related to $R_2$ and $R_3$ given in Lemma~\ref{R2R3}. Here $w=(e-1)q-p$.}
\begin{center}
\begin{tabular}{lllll}
$p$&$q$&$R_2$&$R_3$&$\Delta$\\
\hline\\[-10 pt]
odd&odd&$1^\alpha=3$&$1^\beta=3$&$\Delta_{\bf 01}\alpha=\left\{\begin{array}{ll}\langle w, 2q-2\rangle&e \text{ even}\\
\langle w, -2q\rangle&e \text{ odd}\end{array}\right.$\\[10 pt]
&&&&$\Delta_{\bf 01}\beta=\left\{\begin{array}{ll}\langle w, -2q\rangle&e \text{ even}\\
\langle w, 2q-2\rangle&e \text{ odd}\end{array}\right.$\\[10 pt]
\hline\\[-10 pt]

even&odd&$1^\alpha=3$&$1^\beta=3$&$\Delta_{\bf 01}\alpha=\left\{\begin{array}{ll}\langle w, -2q\rangle&e \text{ even}\\
\langle w, 2q-2\rangle&e \text{ odd}\end{array}\right.$\\[10 pt]
&&&&$\Delta_{\bf 01}\beta=\left\{\begin{array}{ll}\langle w, 2q-2\rangle&e \text{ even}\\
\langle w, -2q\rangle&e \text{ odd}\end{array}\right.$\\[10 pt]
\hline\\[1 pt]

odd&even&$3^\alpha=3$&$1^\beta=1$&$\Delta_{\bf 01}\alpha=\langle w, 2q\rangle$\\[5 pt]
&&&&$\Delta_{\bf 01}\beta=\left\{\begin{array}{ll}\langle w, 2q-2\rangle&e \text{ even}\\
\langle w, -2q-2\rangle&e \text{ odd}\end{array}\right.$\\[10 pt]
\hline

\end{tabular}
\end{center}
\label{R2R3 displacements}
\end{table}%

As already mentioned, we plan to change the presentation of $Q_2(L)$. Finding an equivalent presentation for $Q_2(L)$ will depend heavily on the following lemma.

\begin{lemma} \label{path homotopy} If  $X$ and $Y$ are words in the generators $\{1,2,3\}$ having the same displacement, then   the relation $R_1$ implies that $a^X=a^Y$  for all $a\in \widetilde Q_2$.
\end{lemma}

{\bf Proof:}  Because $2^1=2^3$, it follows that $x^{2^1}=x^{2^3}$ or equivalently that $x^{121}=x^{323}$ for any element $x$ in $\widetilde Q_2$. More generally, suppose $s$ and $t$ are words in the generators 1, 2, and 3 such that their concatenation $st$ is a cyclic permutation of 121323 or its reverse. Then,  for all $x$ and all words $\alpha$ and $\beta$, the relations $x^{\alpha s \beta}=y$ 
and $x^{\alpha {\overline t}\beta}=y$ are equivalent. This last statement can be interpreted in terms of the graph $ {\cal T}_i$ by seeing that we are moving the path $\alpha s \beta$ across one of the hexagonal faces to reach the path  $\alpha {\overline t}\beta$. Call such a move an {\it hexagonal move}. Because $X$ and $Y$ have the same displacement, the paths  starting at $x$ that they each determine clearly end at the same vertex $y$. Moreover, we may move one path to the other by a sequence of hexagonal moves. \hfill $\square$

It is not difficult to show that $A1A232B1B$ and $B1B232B1B$ have the same displacement. Thus, Lemma~\ref{simplify relation}, the proof of which we had postponed,   now follows easily from Lemma~\ref{path homotopy}.

We are now in a position to derive a new presentation of $Q_2(L)$ using Lemma~\ref{path homotopy}. The new presentation will be defined in terms of the following sets of relations. Recall that  $w=(e-1)q-p$ and notice that $R_1$ is equivalent to $\iota_0$.
\begin{align*}
\gamma: 1&=1^{[23]^{|w|}[13]^{(q-2)/2}}&	
\gamma': 3&=3^{[12]^{|w|}[31]^{(q-2)/2}}\\
\delta_i: 1&=1^{[21]^i [31]^q[12]^i}& 				
\delta': 3&=3^{[13]^q} \\
\epsilon_i: 1&=1^{[23]^i[13]^{\pm q/2}[23]^{|w|-i}3}&					
\epsilon_i': 3&=3^{[12]^i[13]^{\pm q/2}[12]^{|w|-i}1}\\
\zeta_{ik}: 1&=1^{[23]^i[31]^k21[31]^{k-1}[32]^{i+1}}&
\zeta_{ik}': 3&=3^{[12]^i[31]^k32[31]^{k+1}[21]^{i+1}}\\
\eta_k: 1&=1^{[23]^{(|w|-1)/2}[31]^k2[13]^{k+q/2}[32]^{(|w|-1)/2}}&
\eta'_k: 3&=3^{[12]^{(|w|-1)/2}[13]^k121[31]^{k+q/2}[21]^{(|w|-1)/2}}\\	
\theta_{ik}: 1&=1^{[21]^i[31]^k21[31]^k[12]^{i+1}}&
\mu: 3&=1^{[21]^{|w|}[31]^{(q-1)/2}}\\
\iota_i: 2&=2^{[12]^i31[21]^i}&
\kappa: 2&=2^{[12]^{|2w|}}\\
\end{align*}

\vspace{-.75 in}
\begin{lemma}\label{new presentation} Let $w=(e-1)q-p$. Then
$$Q_2(L)\cong \left\{ \begin{array}{ll}
\langle 1, 2, 3\,|\, \delta_0, \delta', \gamma, \gamma', \{\epsilon_i, \epsilon_i'\}_{i \in I_2}, \{\zeta_{ik}, \zeta_{ik}'\}_{i\in I_3, k\in I_4}, \{\eta_k, \eta_k'\}_{k \in I_5}, \{\iota_i\}_{i \in I_1}, \kappa \rangle&\text{ $q$ even,}\\
\\
\langle 1, 2, 3\,|\,\{\delta_i\}_{i \in I_1}, \delta', \{\theta_{ik}\}_{i \in I_1, \, k \in I_4 },  \{\iota_i\}_{i \in I_1}, \kappa,\mu   \rangle&\text{ $q$ odd,}
\end{array} \right.
$$
where  $I_1=\{i\,|\,0\le i<|w|\}$, $I_2=\{i\,|\,0< i\le (|w|-1)/2\}$, $I_3=\{i\,|\,0\le i<(|w|-1)/2\}$, $I_4=\{i\,|\,0\le i<q\}$, and $I_5=\{i\,|\,0\le i<q/2\}$.

\end{lemma}
{\bf Proof:} 
There are several cases to consider.  Assume first that $q$ is odd and $w >0$.   We first show that that $R_1,  \delta'$ and $\mu$ can be derived from $R_1, R_2$, and $R_3$ and vice versa. From Table~\ref{R2R3 displacements}, we have the relations $R_2$ and $R_3$ are $1^\alpha=3$ and $1^\beta=3$, respectively. Regardless of the parity of $p$ and $e$, using Lemma~\ref{displacements} one can show that the displacement of $[13]^q$ is the same as that of $\overline{\beta} \alpha 3$ and that the displacement of $[21]^{|w|}[31]^{(q-1)/2}$ is the same as that of $\alpha$.  Therefore, by Lemma~\ref{path homotopy} we have $3^{[13]^q}=3^{\overline{\beta}\alpha 3}$ and $1^{[21]^{|w|}[31]^{(q-1)/2}}=1^\alpha$. Note that these two relations only depend on relation $R_1$.  Therefore, we have:
$$
\begin{array}{lcl}
3 \stackrel{R_2}{=} 1^\alpha \stackrel{R_1}{=} 1^{[21]^{|w|} [31]^{(q-1)/2}} & & \mbox{$\mu$ is a consequence of $R_1$ and $R_2$} \\
3=3^3 \stackrel{R_2}{=}1^{\alpha 3} \stackrel{R_3}{=} 3^{\overline{\beta}\alpha 3} \stackrel{R_1}{=}3^{[13]^q} & & \mbox{$\delta'$ is a consequence of  $R_1, R_2$, and $R_3$} \\
1^\alpha \stackrel{R_1}{=} 1^{[21]^{|w|} [31]^{(q-1)/2}} \stackrel{\mu}{=}3 & & \mbox{$R_2$ is a consequence of $R_1$ and $\mu$} \\
1^\beta \stackrel{R_2}{=} 3^{\overline{\alpha}\beta}=3^{3\overline{\alpha}\beta} \stackrel{\delta'}{=} 3^{[13]^q3\overline{\alpha}\beta} \stackrel{R_1}{=}3 & & \mbox{$R_3$ is a consequence of $R_1$, $\mu$, and $\delta'$}
\end{array}
$$

Assume next that $q$ is odd and $w<0$.  In this case, the displacement of $[13]^q$ is still the same as that of $\overline{\beta} \alpha 3$ but the displacement of $[21]^{|w|}[31]^{(q-1)/2}$ is now the same as that of $1 \beta$.  So, in this case, we use the consequences $3^{[13]^q}=3^{\overline{\beta}\alpha 3}$ and $1^{[21]^{|w|}[31]^{(q-1)/2}}=1^\beta$ of $R_1$.
$$
\begin{array}{lcl}
3 \stackrel{R_3}{=} 1^\beta \stackrel{R_1}{=} 1^{[21]^{|w|} [31]^{(q-1)/2}} & & \mbox{$\mu$ is a consequence of $R_1$ and $R_3$} \\
3=3^3 \stackrel{R_2}{=}1^{\alpha 3} \stackrel{R_3}{=} 3^{\overline{\beta}\alpha 3} \stackrel{R_1}{=}3^{[13]^q} & & \mbox{$\delta'$ is a consequence of  $R_1, R_2$, and $R_3$} \\
1^\beta \stackrel{R_1}{=} 1^{[21]^{|w|} [31]^{(q-1)/2}} \stackrel{\mu}{=}3 & & \mbox{$R_3$ is a consequence of $R_1$ and $\mu$} \\
1^\alpha \stackrel{R_3}{=} 3^{\overline{\beta}\alpha} \stackrel{R_1}{=} 3^{[13]^q 3} \stackrel{\delta'}{=} 3^3 = 3 & & \mbox{$R_2$ is a consequence of $R_1$, $\mu$, and $\delta'$}
\end{array}
$$

Next we will show that $\{\delta_i\}_{i \in I},\{\theta_{i,j}\}_{i \in I, \, j \in J }, \{\iota_i\}_{i \in I}$, and $ \kappa$ can all be derived from $R_1, \mu$, and $ \delta'$. Because $R_2$ and $R_3$ can be derived from $R_1, \mu$, and $ \delta'$, we are free to use these relations as well.

To derive the relations $\{\delta_i\}$, first note that $[21]^i[31]^q[12]^i$ and $\alpha 3 \overline{\beta}$ have the same displacement if $i$ is even. Thus $1^{[21]^i[31]^q[12]^i}=1^{\alpha 3 \overline{\beta}}
=3^{3 \overline{\beta}}=3^{\overline{\beta}}=1$. 
If $i$ is odd, then $[21]^i[31]^q[12]^i$ and $\beta 3 \overline{\alpha}$ have the same displacement and a similar argument shows that $\delta_i$ is a consequence of $R_1, R_2$, and $R_3$.

It is not hard to show that for all $i \in I$ and all $j \in J$, the displacement of $[21]^i [31]^j 21[31]^j[12]^{i+1}$ is $\langle0, 0 \rangle$. Thus the relations $\theta_{i,j}$ all follow from $R_1$.

The word $[12]^i31[21]^i$ has the same displacement as $31$ if $i$ is even and as $13$ if $i$ is odd. Thus, when $i$ is even, we have $2=2^{31}=2^{[12]^i31[21]^i}$, and when $i$ is odd, $2=2^{13}=2^{[12]^i31[21]^i}$. Therefore the relations $\iota_i$ are all consequences of $R_1$. A similar argument actually shows that if $2^X$ is any vertex in ${\cal T}_2$, then $2^{X13}=2^X$ in $\widetilde{Q}_2$. To see this, note that  for all $X$, $\overline{X}31X$ has the same displacement as either $13$ or $31$. So, either 
$$2^{X13}=2^{X\overline{X}31X}=2^{31X}=2^X \quad \text{ or } \quad
2^{X31}=2^{X\overline{X}31X}=2^{31X}=2^X.$$
Thus any two vertices in ${\cal T}_2$ separated by a vertical translation of $4l$ units for any integer $l$ represent equal elements in $\widetilde{Q}_2$.

 Finally, we can easily compute that  $\Delta_{\bf 01}[12]^{|2w|}=\langle -|2w|, 0 \rangle$.
Moreover, the primary relation $1^\alpha=3$ implies that $2^{\overline{\alpha}1 \alpha 3}=2$. Using Lemma~\ref{displacements}, it follows that $\Delta_{\bf 01}\overline{\alpha}1 \alpha 3=-\langle (-1)^w2w, 4q \rangle$. Hence, letting $X$ be either $\overline{\alpha}1 \alpha3$ or its reverse $3\overline{\alpha}1\alpha$, we may assume $\Delta_{\bf 01}X=\langle -|2w|, \pm 4q \rangle$ and hence $\Delta_{\bf 01}X[13]^{\pm q}=\langle -|2w|, 0 \rangle$. Therefore, $2^{[12]^{|2w|}}=2^{X[13]^{\pm q}}=2^{[13]^{\pm q}}=2$ because $[13]^{q}$ represents vertical translation by a multiple of 4. Thus $\kappa$ is a consequence of $R_1$ and $R_2$.

This proves the Lemma in the case where $q$ is odd. We leave the proof in the case where $q$ is even to the reader. \hfill $\square$

{\bf Proof of Theorem~\ref{Cayley graph}:} We will now apply Winker's method to the  presentation in Lemma~\ref{new presentation}, focussing first on the case where $q$ is odd.   As we follow the method, imagine filling in vertices and edges in the graphs ${\cal C}_1$ and ${\cal C}_2$ already described in Figures~\ref{C1 q odd} and \ref{C2}. We start with the introduction of three vertices representing the generators 1, 2, and 3, and each with a loop attached due to the idempotentcy relations.  Next, the primary relations are used to add additional vertices and edges, if necessary. Moreover, we apply these in the following order: $\mu, \delta', \delta_i, \theta_{ik}, \kappa$, and $\iota_i$. Tracing $\mu$ connects the generator 1 to the generator 3 along a horizontal section (a path of edges $[21]^{|w|}$) followed by the right vertical side.  Next, tracing $\delta'$ closes the loop at the lower right corner. The  relation $\delta_0$ fills in the left vertical side. The relations $\delta_i$, for $i>0$,  now leave the horizontal section, fill in the vertical sections, and then return to the horizontal section. Finally,  the relations $\theta_{ik}$  fill in all the remaining horizontal edges labeled 2. The component ${\cal C}_1={\cal C}_3$ is now complete. The relation $\kappa$ traces out the horizontal section of ${\cal C}_2$ and finally, the relations $\iota_i$ fill in the vertical sections, thus completing ${\cal C}_2$.

The next step in Winker's method is to apply the secondary relations to each already existing vertex. In general, when  checking a secondary relation for a given vertex $x$, it may be necessary to introduce more vertices and edges or cause collapsing to take place. However, in our case, we will see  that each secondary relation is already satisfied at each vertex. We chose this presentation of $Q_2(L)$ precisely for this reason. 

The secondary relations, which must be true for all $x$,  are listed below. We label each with the associated name of the primary relation preceded by ``s'', for secondary. The indices range over the same index sets as the corresponding primary relation.
\begin{align*}
s\mu: x&=x^{[13]^{(q-1)/2}[12]^{|2w|}[13]^{(q+1)/2}}\\
s\delta': x&=x^{[31]^{2q}}\\
s\delta_i: x&=x^{[21]^i[13]^q[12]^{2i}[13]^q[21]^i}\\
s\theta_{ik}: x&=x^{[21]^{i+1}[13]^k12[13]^k[12]^{2i}[13]^k12[13]^k[21]^{i+1}}\\
s\kappa: x&=x^{[21]^{|4w|}}\\
s\iota_i: x&=x^{[12]^i13[21]^{2i}23[12]^{i+1}}
\end{align*}
To examine the effect of the secondary relations at every element $x$ we compute the following displacements:
\begin{align*}
\Delta_{\bf01}[13]^{(q-1)/2}[12]^{|2w|}[13]^{(q+1)/2}&=-\langle 2|w|, 4q\rangle\\
\Delta_{\bf01}{[31]^{2q}}&=\langle0, 8q \rangle\\
\Delta_{\bf01}{[21]^i[13]^q[12]^{2i}[13]^q[21]^i}&=(-1)^i\langle0, -8q \rangle\\
\Delta_{\bf01}[21]^{i+1}[13]^k12[13]^k[12]^{2i}[13]^k12[13]^k[21]^{i+1}&=\langle0,0 \rangle\\
\Delta_{\bf01}[21]^{|4w|}&=\langle4|w|,0 \rangle\\
\Delta_{\bf01}[12]^i13[21]^{2i}23[12]^{i+1}&=\langle0,0 \rangle\\
\end{align*}
Because horizontal translation by $2|w|$ and vertical translation by $4q$ preserve elements of the quandle in both ${\cal T}_1$ and ${\cal T}_2$, it follows that the secondary relations hold at every vertex in ${\cal C}_1$ and ${\cal C}_2$.

This completes the proof of Theorem~\ref{Cayley graph} in the case where $q$ is odd. We leave the case when $q$ is even to the reader. In working through the details, notice that the construction of ${\cal C}_1$ and  ${\cal C}_3$ run in parallel fashion and that the construction of ${\cal C}_2$ is the same as in the case where $q$ is odd.
\hfill $\square$

\pagebreak
\section{Some Properties of $Q_2(L)$}\label{properties}

Two-bridge links are not distinguished by their involutory quandles because $Q_2(K_{p/q})\cong \text{Core}(\mathbb Z/q \mathbb Z)$, see \cite{JO2}. However, for $(2, 2, r)$-Montesinos links we have the following theorem.

\begin{theorem} If $L$ and $L'$ are both $(2, 2, r)$-Montesinos links with isomorphic involutory quandles, then $L$ is either $L'$ or its mirror image.
\end{theorem}
{\bf Proof:} Let $L=L(1/2,1/2,p/q;e)$ and $L'=L(1/2,1/2,p'/q';e')$. Assume that $Q_2(L')\cong Q_2(L)$,  $w=(e-1)q-p$, and $w'=(e'-1)q'-p'$. Using Corollary~\ref{corollary} to compute  the number of components of the quandles, as well as their sizes we  conclude that $|w|=|w'|$ and $q=q'$. If $(e-1)q-p=(e'-1)q-p'$, then $q$ divides $p-p'$ from which it follows that $p=p'$ and $e=e'$. Hence the links are equivalent. On the other hand, if $(e-1)q-p=-(e'-1)q+p'$, then $q$ divides $p+p'$ and we must have $p+p'=q$. This gives that $e=-e'+3$. But these conditions now imply that $L$ and $L'$ are mirror images of each other. \hfill $\square$

In any involutory quandle $Q$, the subquandle generated by two elements $x$ and $y$ is called the {\it geodesic} through $x$ and $y$, which we denote $g(x,y)$. A {\it maximal} geodesic is one that is not a proper subset of any other geodesic. More details regarding geodesics are given in \cite{JO}, \cite{JO2}. The following theorem describes the set of maximal geodesics in any $(2,2,r)$-Montesinos link.

\begin{theorem} 
Let $L=L(1/2, 1/2, p/q;e)$ and $w=(e-1)q-p$. The involutory quandle $Q_2(L)$ has $q+1$ maximal geodesics. 
\vspace{-.2 in} 
\begin{itemize}
\item[(a)] If $q$ is odd, then  one maximal  geodesic is $g(1, 1^{32})$ of length $2 q|w|$ and  the remaining $q$ maximal geodesics  are $g(2, 3), g(2, 3^1), \dots, g(2, 3^{[13]^{(q-1)/2}})$,  each of length $4|w|$. The geodesic  $g(1, 1^{32})$ is equal to ${\cal C}_1$ and the intersection of any two of the other geodesics is ${\cal C}_2$.
\item[(b)] If $q$ is even, then one  maximal geodesic is $g(1,3^2)$ of length $2q|w|$ and the remaining $q$ maximal geodesics are $g(2, 1)$, $g(2, 3)$, $g(2, 1^3)$, $g(2, 3^1)$, $g(2, 1^{31})\dots$,  each of length $4|w|$. The geodesic  $g(1,3^2)$ is equal to ${\cal C}_1 \cup {\cal C}_3$ and the intersection of any two of the other geodesics is ${\cal C}_2$.
\end{itemize}
\end{theorem}
{\bf Proof:} Suppose $q$ is odd and consider the geodesic generated by $1$ and $1^{32}$. Traveling only ``forward'' on this geodesic, we generate the sequence
$$1, 1^{32}, 1^{123132}, 1^{32123132}, \dots, 1^{[123132]^k}, 1^{32[123132]^k}, \dots.$$
Note that $\sigma(123132)=\text{id}$ and that $\Delta_{\bf 01}[123132]^k=-2k\langle 1, 4 \rangle$. Thus, the element $1^{[123132]^k}$ is represented by the vertex $(-2k, -8k+1)$ in ${\cal T}_1$. This is a vertex of type $\bf 01$. Similarly, the element $1^{32[123132]^k}$ is represented by the vertex $(-2k-1, 8k+3)$, a vertex of type {\bf 13}. 

Now consider moving any  vertex $P$ of ${\cal T}_1$  to an equivalent vertex which lies in the fundamental domain ${\cal F}_1$. We can do this by first translating $P$ to $P'$, a vertex in the domain $[0, 2w)\times (0, 4q)$, by  using translation in the horizontal direction by a multiple of $2w$ and in the vertical direction by a multiple of $4q$. If $P' \not\in {\cal F}_1$, then we may reflect $P'$ through the point $(w, 2q)$ to obtain the point $P'' \in {\cal F}_1$. 

If $k=q|w|$, then the element $1^{[123132]^{q|w|}}$ located at $(-2q|w|, -8q|w|+1)$ is equivalent to the vertex $(0, 1)$ which is the element $1$. Similarly, the element $1^{32[123132]^{q|w|}}$ located at $(-2q|w|-1, 8q|w|+3)$ is equivalent to the vertex $(1, 4q-3)$ which is the element $1^{32}$. Thus the geodesic has returned to the beginning after traversing $2q|w|$ points. We want to show that all $2q|w|$ of these points are distinct in ${\cal F}_1$, and therefore make up all of the component ${\cal C}_1$. If not, then there are two distinct points $P$ and $Q$ on the geodesic when thought of as points in ${\cal T}_1$ that project to the same point in ${\cal F}_1$. One possibility is that $P$ and $Q$ translate to the same point $P'=Q'$ in  $[0, 2w)\times (0, 4q)$. Because such a translation preserves the  type of a vertex, we see that $P$ and $Q$ must both be of type {\bf 01} or {\bf 13}. If this were true, there would exist integers $j$ and $k$ with $0\le j<k<q|w|$ such that $P=(-2j, -8j+1)$ and $Q=(-2k, -8k+1)$ or  $P=(-2j-1, 8j+3)$ and $Q=(-2k-1, 8k+3)$. In either case, $P'=Q'$ implies that $j\equiv k$ (mod $|w|$) and $j\equiv k$ (mod $q$) because $q$ is odd. Thus there exists $n$ and $m$ such that $k=j+n|w|$ and $k=j+mq$. Furthermore, $0\le n<q$ and $0\le m<|w|$. From this we obtain that $mq=n|w|=\pm n[(e-1)q-p]$ and so $q$ divides $np$. Because $p$ and $q$ are co-prime and $n<q$, it follows that $n=0$ and $j=k$. Another possibility is that $P'\ne Q'$ but instead $P'$ and $Q'$ are symmetric with respect to either the point $(0,2q)$ or $(|w|,2q)$. But any pair of such  symmetric points must have types $\bf 01$ and $\bf 03$ or types $\bf 11$ and $\bf 13$. This is impossible because $P$ and $Q$ have types {\bf 01} or {\bf 13} and these types are preserved by the translations. Thus, the geodesic contains $2q|w|$ distinct points, filling out all of ${\cal C}_1$.

A similar argument can be used to show that the geodesic generated by $2$ and $3^{[13]^k}$, where $0\le k \le (q-1)/2$ alternates between ${\cal C}_2$ and ${\cal C}_1$ traversing every point of ${\cal C}_2$ exactly once and every point of ${\cal F}_1$ with $y$-coordinate equal to $2q\pm(4k+1)$ exactly once. This is a total of $2|w|$ points in ${\cal C}_1$ and $2|w|$ points in ${\cal C}_2$. We leave this case, as well as the case of the geodesic generated by $2$ and $3^{[13]^k1}$ to the reader. \hfill $\square$  

If $Q$ is a quandle an {\it automorphism} of $Q$ is a bijection $\phi :Q \to Q$ such that $\phi(x \rhd y)=\phi(x) \rhd \phi(y)$ for all $x$ and $y$ in $Q$. The set of automorphism forms a group, $\text{Aut}(Q)$,  under composition in the obvious way.  any automorphism must take a maximal geodesic in an involutory quandle  to a maximal geodesic. Using this idea we may prove the following theorem.

\begin{theorem} \label{order of automorphism group}The order of $\text{Aut}(Q_2(L))$ is at most $24 w^2 \Phi(4|w|)$ if $q=2$ and $4qw^2\Phi(2 q|w|)$ if $q>2$,  where $\Phi$ is the Euler $\Phi$-function.
\end{theorem}
{\bf Proof:} First note that any automorphism must preserve both maximal geodesics and components.
For simplicity, we consider the case where $q$ is odd. The even case is similar. If $q>2$, there is one   maximal geodesic of length $2 q|w|$ which must be taken to itself and $q-1$ shorter geodesics of length $4|w|$ that must be permuted among themselves. 
Since $1$ and $1^{32}$ generate the long geodesic, we must have that this pair of elements is sent to another pair of generators of this geodesic and knowing where they are sent is sufficient to determine the image of any other element of this geodesic. Because the geodesic contains $2q|w|$ elements, there are $2q|w|$ possible choices for $\phi(1)$. In order for $\phi(1)$ and $\phi(1^{32})$ to generate the geodesic, we need these two points located a distance $d$ apart with $d$ prime to $2q|w|$. Thus there are $\Phi(2q|w|)$ possibilities for $\phi(1^{32})$, where $\Phi$ is the Euler $\Phi$-function. We also know that $2$ must be taken to another point of ${\cal C}_2$. Thus there are $2|w|$ possibilities for its image. Once this is determined, the image of every point is determined because each point lies on a maximal geodesic containing 2. Thus there are at most $4qw^2 \Phi(2q|w|)$ elements in $\text{Aut}(Q_2(L))$. The proof for $q=2$ is similar.
 \hfill $\square$

Experimental evidence suggests that the order of $\text{Aut}(Q_2(L))$ exactly equals the upper bound given in Theorem~\ref{order of automorphism group}.

 \end{document}